\documentclass[11pt]{article}
\usepackage{bez123,calc,curves,ebezier,epic,eepic,graphicx,multiply,rotating}
\usepackage{epic,eepic,graphicx}
\usepackage{tikz}
\everymath{\displaystyle}
\usepackage{graphicx}
\usepackage{graphicx}
\usepackage{pst-all}
\usepackage{color}
\usepackage{multirow}
\usepackage{amssymb}
\usepackage{amsmath}
\newtheorem{prethm}{{\bf Theorem}}

\newenvironment{thm}{\begin{prethm}{\hspace{-0.5
				em}{\bf .}}}{\end{prethm}}
\newtheorem{prelemma}{{\bf Lemma}}

\newtheorem{preex}{{\bf Example}}

\newtheorem{preprop}{{\bf Proposition}}

\newenvironment{prop}{\begin{preprop}{\hspace{-0.5em}{\bf .}}}{\end{preprop}}
\newtheorem{precor}{{\bf Corollary}}

\newenvironment{cor}{\begin{precor}{\hspace{-0.5
				em}{\bf .}}}{\end{precor}}
\newtheorem{preremark}{{\bf Remark}}

\newtheorem{preprob}{{\bf Problem}}

\newtheorem{predefin}{{\bf Definition}}

\newenvironment{defin}{\begin{predefin}{\hspace{-0.5
				em}{\bf .}}}{\end{predefin}}
\newtheorem{preconj}{{\bf Conjecture}}

\newtheorem{preprobb}{{\bf Problem}}

\newtheorem{prelem}{{\bf Theorem}}

\newenvironment{proof}{{\bf Proof.}\rm }{\hfill{$\Box$}}

\newtheorem{presolution}{{\bf Solution.}}

\def\emline#1#2#3#4#5#6{%
\put(#1,#2){\special{em:moveto}}%
\put(#4,#5){\special{em:lineto}}}
\def\newpic#1{}

\title{\vspace{-0.1cm}\Large\bf More relations between $\lambda$-labeling and Hamiltonian paths with emphasis on line graph of bipartite multigraphs}
\author{\large\bf Manouchehr Zaker\footnote{mzaker@iasbs.ac.ir}
\vspace{5mm}\\
Department of Mathematics,\\ Institute for Advanced Studies
in Basic Sciences,\\ Zanjan 45137-66731, Iran
}
\date{}
\begin{document}	
\maketitle
\begin{abstract}
\noindent This paper deals with the $\lambda$-labeling and $L(2,1)$-coloring of simple graphs. A $\lambda$-labeling of a graph $G$ is any labeling of the vertices of $G$ with
different labels such that any two adjacent vertices receive
labels which differ at least two. Also an $L(2,1)$-coloring of $G$ is any labeling of the vertices of $G$ such that any two adjacent vertices receive
labels which differ at least two and any two vertices with distance two receive distinct labels. Assume that a partial $\lambda$-labeling $f$ is given in a graph $G$. A general question is whether $f$ can be extended to a $\lambda$-labeling of $G$. We show that the extension is feasible if and only if a Hamiltonian path consistent with some distance constraints exists in the complement of $G$. Then we consider line graph of bipartite multigraphs and determine the minimum number of labels in $L(2,1)$-coloring and $\lambda$-labeling of these graphs. In fact we obtain easily computable formulas for the path covering number and the maximum path of the complement of these graphs. We obtain a polynomial time algorithm which generates all Hamiltonian paths in the related graphs. A special case is the Cartesian product graph $K_n\Box K_n$ and the generation of $\lambda$-squares.
\end{abstract}



\section{Introduction}

\noindent All graphs in this paper are undirected graphs and without any loops. Let $G=(V(G),E(G)$ be a graph without any multiple edges. A $\lambda$-labeling of $G$ is any labeling of the vertices of $G$ with
different labels such that any two adjacent vertices receive
labels which differ at least two. This concept comes from
assigning non-reusable frequencies to radio transmitters in such a
way that close transmitters use frequency channels which are
sufficiently separated from each other and is widely studied in
literature \cite{BKTV, FKK, FNPS, FS, Y}. Another type of distance constrained labeling for graphs is the so-called $L(p,q)$-labeling, where $p$ and $q$ are non-negative integers. An $L(p,q)$-labeling of $G$ is any assignment of non-negative integers to the vertices of $G$ such that $(i)$ any two adjacent vertices receive labels which differ at least $p$ and $(ii)$ any two vertices with
distance two receive labels which differ at least $q$. In this paper instead of the term $L(2,1)$-labeling we use the term $L(2,1)$-coloring which is also used in many papers. The $L(2,1)$-coloring is also called radio coloring in some papers and widely studied in the literature. It is clear that for graphs having diameter two any $\lambda$-labeling is also an $L(2,1)$-coloring and vice versa. The minimum $\lambda$ for which
there exists an $L(2,1)$-coloring of $G$ using labels from $\{0, 1, \ldots, \lambda\}$ is denoted by
$\lambda_{2,1}(G)$. The literature is full of papers concerning $L(p,q)$-labeling and in particular $L(2,1)$-coloring of graphs \cite{FNPS, GMS, GY, HCHY, KY, LLS, SS, SYZ, Y}. Algorithmic results on radio colorings have been obtained by Bodlaender et al. in \cite{BKTV}. In this paper we also consider the concept of partial $\lambda$-labeling.
Let $S$ be any subset of vertices in $G$ and $f$ be an assignment of distinct non-negative integers to the vertices of $S$ such that if $u, v$ are two adjacent vertices of $S$ then $|f(u)-f(v)|\geq 2$. Then the subset $S$ together with $f$ is called a partial $\lambda$-labeling of $G$. The question whether a partial $\lambda$-labeling (resp. partial radio labeling) can be extended to a $\lambda$-labeling (resp. radio coloring) of whole graph has been studied in many papers (e.g. \cite{BBFPW,FKK,FKP1}). In \cite{FKP1} a polynomial time algorithm for the extension problem of partial $\lambda$-labelings of trees has been obtained. In this paper we use the concept of line graphs. Let $G$ be
any graph with or without multiple edges. Let $E(G)$ be the edge set of $G$. The line graph $L(G)$ of $G$ is a graph whose vertex set is $E(G)$ in which any two edges $e_1, e_2 \in E(G)$ (as two vertices in $L(G)$) are adjacent
if and only if they have at least one common end-vertex.

\noindent The other subject to be studied in the paper is Hamiltonian path, because there are close relationships between $\lambda$-labeling and Hamiltonian paths. In a graph $G$ with $n$ vertices any path consisting of $n$ vertices is called a Hamiltonian path. This relationship is appeared first time in
\cite{GY}, where the following theorem is proved. In this paper the complement of any graph $G$ is denoted by $\overline{G}$.

\begin{thm}\cite{GY}\label{first}
In any graph $G$ on $n$ vertices the following two statements are
equivalent:

\noindent 1. There exists an injective $f: V(G) \rightarrow
\{0,1,\ldots,n-1\}$ such that $|f(x)-f(y)|\geq 2$ for all
$\{x,y\}\in E(G)$;\\
\noindent 2. $\overline{G}$ contains a Hamilton path.
\end{thm}

\noindent The path covering number $\pi(G)$ of a graph $G$ is the smallest number of vertex disjoint paths needed to cover the vertices $G$. A result of Georges et al. \cite{GMW} relates path covering number of $\overline{G}$ to $\lambda_{2,1}(G)$.

\begin{thm}\cite{GMW}\label{path cov}
Let $G$ be a graph on $n$ vertices with $\pi(\overline{G})=r\geq 2$. Then
$$\lambda_{2,1}(G)=n+r-2.$$
\end{thm}

\noindent The study of relations between Hamiltonian paths and path covering number of graphs
in one side and distance constrained labeling of graphs such as radio colorings in the other side, has been the research subject of many articles \cite{FS, GMW, LZ, SA}.

\noindent {\bf The paper is organized as follow.} In Section $2$ we consider partial $\lambda$-labeling of graphs and study the problem whether a partial $\lambda$-labeling can be extended to a $\lambda$-labeling of the entire graph. We relate this problem to a new version of Hamiltonian path problem (Theorem \ref{equi}). In Section $3$ we study the $\lambda$-labeling (or $L(2,1)$-coloring) of line graph of bipartite multigraphs. The most important notions that are used are pathwise tough graph and tough vertex. For simplicity the complement of a line graph of any bipartite multigraph is called skew graph. In Section $3$ pathwise tough skew graphs and tough vertices have been characterized (Theorem \ref{skew-pathwise} and Corollary \ref{tough}). In Section $4$ we prove that a connected skew graph $G$ admits a Hamiltonian path from a vertex $v$ if and only if $v$ is tough in $G$ (Theorem \ref{hamilton}). Also a skew graph contains a Hamiltonian path if and only if it is pathwise tough (Corollary \ref{hamilton-pathwise}). An efficient algorithm is given which obtains (if exists) and generates all Hamiltonian paths in skew graphs (Corollary \ref{generate}). Also the path covering number and the length of largest path of skew graphs are determined by easy formulas (Proposition \ref{pathcov} and Theorem \ref{maxpath}).

\section{Hamiltonian paths consistent with a pre-labeling of vertices}

\noindent In this section we study extendibility of partial $\lambda$-labellings using concepts related to Hamiltonian paths in graphs. In the following we introduce two concepts which are closely related. First, we need to recall some notation. Let $H$ be a graph and $u, v$ are two vertices in $H$. Denote by $d_H(u,v)$ the length of smallest path in $H$ between $u$ and $v$.

\begin{defin}
Let $G$ be a graph and $S=\{v_1, \ldots, v_k\}$ any (ordered) subset of vertices in $G$. By a distance pattern ${\mathcal{D}}$ for $S$ we mean any function on $\{(v_1, v_2), \ldots, (v_i, v_{i+1}), (v_{k-1}, v_k)\}$ which for each $i=1, \ldots, k-1$ assigns a positive integer $d_i$ to the pair $(v_i, v_{i+1})$. Let $P$ be any Hamiltonian path in $G$. We say that $P$ is consistent with the distance pattern ${\mathcal{D}}$ if the distance of $v_{i}$ and $v_{i+1}$ (along the path $P$) equals to $d(v_{i},v_{i+1})$. In other words for each $i\in \{1, \ldots, k-1\}$, $d_P(v_i,v_{i+1}) =d_i.$\label{pathformul1}
\end{defin}

\noindent This subject is clearly closely related to the disjoint
path problem with prescribed distance between the given pairs of
vertices. In the ordinary disjoint paths problem, we are given a graph
$G$ and a set consisting of some pairs of vertices $\{(s_i,t_i)\}$, $i=1, \ldots, k$. The question is whether
there are $k$ vertex disjoint paths connecting $s_i$ to $t_i$ for each $i$. In
its generalized form we seek for $k$ vertex disjoint paths $P_1, \ldots, P_k$ such that $P_i$ is between $s_i$ and $t_i$ and its length is a prescribed number $d_i$.

\noindent Definition \ref{pathformul1} has another formulation which is used in the paper.
Given any set of $n$ cities and some roads between them, let a graph $G$ represent the cities and roads in which a road between two arbitrary cities $c_i$ and $c_j$ is represented by an edge between the vertices $c_i$ and $c_j$ in $G$. Let $S$ be a subset of cities in $G$. Let $\tau: S\rightarrow
\{0, \ldots, n-1\}$ be an assignment of labels to the elements of $S$. A tourist plans to start from a city $u$ and move along the roads such that each city $v\in S$ is visited in the $\tau(v)$-th round (or day) of her trip. It is assumed that to pass from one edge (road) takes one round or day in the trip. Here, we don't need to consider costed edges but if we associate a cost to each edge then a new variant of Traveling Salesman Problem is emerged in which the salesman plans to visit some cities $c_i$ in $t_i$-th round of her trip, where $t_i$ is a preassigned integer associated to $c_i$. This formulation of Definition \ref{pathformul1} leads us to introduce the following concept. Given a graph $G$, a subset of cities $S$ and the function $\tau$ as above, we say a Hamiltonian path $P$ in $G$ with
starting point $u$ is {\it consistent} with $(S,\tau)$
if the distance of any vertex $v\in S$ from $u$ in the path $P$ is
$\tau(v)$. In other words, if we start from $u$ and move along
the path $P$ then we visit any vertex $v$ from $S$ in $\tau(v)$-th
step.

\begin{defin}
\noindent Let $G$ be a graph and $S$ a subset of vertices in $G$. Let $\tau: S\rightarrow \{0, \ldots, n-1\}$ be an assignment of distinct labels to the vertices of $S$. A $u$-path $P$ in
$G$ is said to be a $\tau$-consistent path, if the distance of
$v$ from $u$ in the path $P$ is $\tau(v)$, for each $v\in S$.
\end{defin}

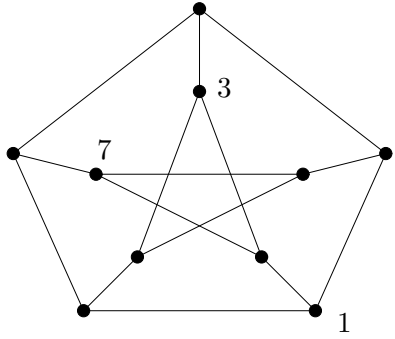
\begin{figure}
\centering \unitlength=0.55mm \special{em:linewidth 0.4pt}
\linethickness{0.7pt}
\begin{picture}(65.00,84.17)
\put(7.00,7.00){\circle*{3.00}}
\put(63.00,7.00){\circle*{3.00}}
\put(-10.00,45.00){\circle*{3.00}}
\put(80.00,45.00){\circle*{3.00}}
\put(35.00,80){\circle*{3.00}}
\put(35.00,60.00){\circle*{3.00}}
\put(10.00,40.00){\circle*{3.00}}
\put(60.00,40.00){\circle*{3.00}}
\put(20.00,20.00){\circle*{3.00}}
\put(50.00,20.00){\circle*{3.00}}
\emline{7.00}{7.00}{1}{63.00}{7.00}{2}
\emline{63.00}{7.00}{1}{80.00}{45.00}{2}
\emline{80.00}{45.00}{1}{35.00}{80.00}{2}
\emline{35.00}{80.00}{1}{-10.00}{45.00}{2}
\emline{-10.00}{45.00}{1}{7.00}{7.00}{2}
\emline{20.00}{20.00}{1}{35.00}{60.00}{1}
\emline{50.00}{20.00}{1}{10.00}{40.00}{1}
\emline{50.00}{20.00}{1}{35.00}{60.00}{2}
\emline{60.00}{40.00}{2}{10.00}{40.00}{2}
\emline{60.00}{40.00}{1}{20.00}{20.00}{2}
\emline{7.00}{7.00}{1}{20.00}{20.00}{2}
\emline{63.00}{7.00}{1}{50.00}{20.00}{2}
\emline{80.00}{45.00}{1}{60.00}{40.00}{2}
\emline{35.00}{80.00}{1}{35.00}{60.00}{2}
\emline{-10.00}{45.00}{1}{10.00}{40.00}{2}
\put(70.00,4.00){\makebox(0,0)[cc]{$1$}}
\put(41.00,61.00){\makebox(0,0)[cc]{$3$}}
\put(12.00,46.00){\makebox(0,0)[cc]{$7$}}
\end{picture}
\caption{Does there exist a Hamiltonian path in the Petersen graph
consistent with the above labels?}\label{peter3}
\end{figure}

\noindent For example in Figure \ref{peter3}, three vertices with
labels are given. There are two paths of length four between
vertices $3$ and $7$. But no Hamiltonian path is consistent with
these labels.

\noindent The following theorem shows a relationship between newly
defined concepts and the former ones. In the proof we use the concept of internally vertex disjoint paths. Given a path graph $P$, a vertex $w$ is called an end-vertex of $P$ is the degree of $w$ in $P$ is one. In a graph $G$, a collection of paths $P_1, P_2, \ldots, P_k$ is said to be internally vertex disjoint if for each $i$ and $j$, $i\not= j$, either $V(P_i)\cap V(P_j)$ is empty set or any vertex in $V(P_i)\cap V(P_j)$ is a common end-vertex of $P_i$ and $P_j$.

\begin{thm}\label{equi}
Let $G$ be a graph on $n$ vertices, $S=\{v_1, \ldots, v_t\}$ a subset of vertices in $G$ and $c:S \rightarrow \{0, 1, \ldots, n-1\}$ a partial $\lambda$-labeling (or radio coloring) of $G$. Assume that $c(v_1)< c(v_2) < \cdots < c(v_t)$. The following statements are equivalent:

\noindent (i) $c$ can be extended to a $\lambda$-labeling (or radio coloring)
of $G$ with $n$ labels.

\noindent (ii) There is a Hamiltonian path in $\overline{G}$ consistent with
$c$.

\noindent (iii) There exist vertex disjoint paths in $\overline{G}$ between
any $v_i$ and $v_{i+1}$ with length $c(v_{i+1}) - c(v_i)$, for
each $i=1,\ldots,t-1$ and a $v_1$-path with length $c(v_1)$ and
a $v_t$-path with length $n-c(v_t)-1$.
\end{thm}

\noindent \begin{proof}
Assume first that $(i)$ holds. Let $f$ be a $\lambda$-labeling of $G$ such that for each $i$, $f(v_i)=c(v_i)$. Let $c(v_i)=c_i$. If $u$ and $v$ are two vertices in $G$ with $|f(u)-f(v)|=1$ then $u$ and $v$ are adjacent in $\overline{G}$. Let $u_0, u_1, \ldots, u_{c_1}$ be such that $f(u_i)=i$, for each $i=0, \ldots, c_1$. We have $u_{c_1}=v_1$ since $f(v_1)=c_1$ and $f$ is injective. It follows that $u_i$ is adjacent to $u_{i+1}$ for each $i$ with $0\leq i \leq c_1-1$. In other words, $\{v\in V(G): 0 \leq f(v) \leq c_1\}$ forms a path in $\overline{G}$ which begins from $u_0$ in $G$ and ends at $v_1$. Note that if $c_1=0$ then $u_0=v_1$. Also the following sets $\{v\in V(G): c_i\leq f(v) \leq c_{i+1}\}$, where $i=1, \ldots, t-1$, are internally vertex disjoint paths in $\overline{G}$, because similar to the case $\{v\in V(G): 0 \leq f(v) \leq c_1\}$ each $\{v\in V(G): c_i\leq f(v) \leq c_{i+1}\}$ is a path and since $f$ is injective then these paths are internally vertex disjoint. Finally, the set $\{v\in V(G): c_t \leq f(v) \leq n-1\}$ forms a path in $\overline{G}$ which begins from $v_t$ and ends at some vertex say $w$ in $G$. Clearly all of these paths are internally disjoint paths in $\overline{G}$, because $f$ is a $\lambda$-labeling. Combining these paths we obtain a Hamiltonian path say $P$ in $\overline{G}$ which starts from $u$ and the distance of $v_i$ from $u$ in $P$ is $c(v_i)=c_i$ for each $i=1, \ldots, t$. This argument shows that $(i)\Rightarrow (ii)$. It can be easily observed that $(ii)$ and $(iii)$ are equivalent. To complete the proof we show that $(ii)$ implies $(i)$. Let $P$ be a Hamiltonian path in $\overline{G}$ consistent with $c$. Let $u$ be the starting vertex of $P$. Define a labeling $f$ for the vertices of $G$ as follows. For each $v_i$, define $f(v_i)=c_i$. In general for any vertex $v$ set $f(v)=d_P(u,v)$. Assume that $v$ and $w$ are two arbitrary vertices in $G$ such that $|f(v)-f(w)|=1$. It follows that $v$ and $w$ are two consecutive vertices in $P$, i.e. they are adjacent in $\overline{G}$. Hence $v$ and $w$ are not adjacent in $G$.
\end{proof}

\section{Skew graphs and Hamiltonian paths: part I}

\noindent In order to begin this section and introduce the concept of skew graphs we need some introduction. Let $G$ be a bipartite graph without any loops and multiple edges, with the bipartite sets $X=\{x_1, \ldots, x_m\}$ and $Y=\{y_1, \ldots, y_n\}$. Note that the line graph $L(G)$ of $G$ can be represented by an $m\times n$ table $T$ of boxes, whose rows and columns are indexed by $1, 2, \ldots, m$ and $1, 2, \ldots, n$, respectively. For any $i, j$ with $1\leq i \leq m$ and $1\leq j \leq n$, if there exists an edge between $x_i$ and $y_j$ in $G$ then we put a vertex in the box corresponding to $(i,j)$ in the table. By the definition of line graphs, any two vertices in $T$ are adjacent in $L(G)$ if and only if they are in a same row or column. Clearly, the vertex set of the complement $\overline{L(G)}$ of $L(G)$ can be represented by the same table $T$ but any two vertices in $T$ are adjacent in $\overline{L(G)}$ if and only if they are in different rows and also in different columns. Now, if we go from simple bipartite graphs $G$ to bipartite graphs $H$ with multiple edges and the same bipartite sets $X$ and $Y$, same tabular representations are obtained for $L(H)$ and $\overline{L(H)}$, only by a slight modifications, i.e. in any $(i,j)$-box we put $t_{ij}$ vertices if there are $t_{ij}$ edges between $x_i$ and $y_j$ in $H$. Any two vertices which are in
a same box are adjacent in $L(H)$ and hence are not adjacent in $\overline{L(H)}$. The concept of skew graphs are obtained by considering this representation for $L(H)$. The details are as follow.

\begin{figure}
\unitlength=0.55mm \special{em:linewidth 0.4pt}
\linethickness{0.7pt}
\hspace*{3.2cm}
\vspace*{0cm}\begin{picture}(75.00,84.17)
\put(25.00,70.00){\circle*{3.00}}
\put(-10.00,45.00){\circle*{3.00}}
\put(30.00,45.00){\circle*{3.00}}
\put(40.00,45.00){\circle*{3.00}}
\put(50.00,45.00){\circle*{3.00}}
\put(60.00,45.00){\circle*{3.00}}
\put(10.00,20.00){\circle*{3.00}}
\put(40.00,20.00){\circle*{3.00}}
\emline{25.00}{70.00}{1}{-10.00}{45.00}{2}
\emline{25.00}{70.00}{1}{30.00}{45.00}{2}
\emline{25.00}{70.00}{1}{40.00}{45.00}{2}
\emline{25.00}{70.00}{1}{50.00}{45.00}{2}
\emline{25.00}{70.00}{1}{60.00}{45.00}{2}
\emline{25.00}{70.00}{1}{10.00}{20.00}{2}
\emline{30.00}{45.00}{1}{10.00}{20.00}{2}
\emline{40.00}{45.00}{1}{10.00}{20.00}{2}
\emline{-10.00}{45.00}{1}{10.00}{20.00}{2}
\emline{40.00}{20.00}{1}{10.00}{20.00}{2}
\emline{50.00}{45.00}{1}{40.00}{20.00}{2}
\emline{60.00}{45.00}{1}{40.00}{20.00}{2}
\emline{-10.00}{45.00}{1}{40.00}{20.00}{2}
\end{picture}
\setlength{\unitlength}{1.8mm}
\hspace*{-5.5cm}
\begin{picture}(90,40)
\put(40,8){\line(1,0){16}}
\put(40,12){\line(1,0){16}}
\put(40,16){\line(1,0){16}}
\put(40,20){\line(1,0){16}}
\put(41,17.5){$1$}
\put(45,17.5){$1$}
\put(45,13.5){$2$}
\put(49,13.5){$2$}
\put(53,13.5){$1$}
\put(49,9.5){$1$}
\large
\put(40,8){\line(0,1){12}}
\put(44,8){\line(0,1){12}}
\put(48,8){\line(0,1){12}}
\put(52,8){\line(0,1){12}}
\put(56,8){\line(0,1){12}}
\small\put(35,2){} \large
\end{picture}
\vspace*{-1cm}
\caption{An example of skew graph and its corresponding table}
\label{fig-skew}
\end{figure}
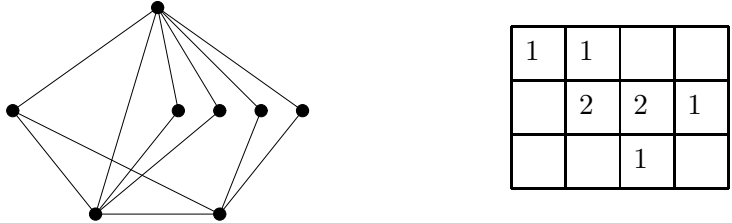

\noindent Let $H$ be any bipartite graph without any loops but possibly with multiple edges (i.e. parallel edges). Let also $L(H)$ be its line graph. In the rest of the paper we study the $\lambda$-labeling and $L(2,1)$-coloring of such graphs $L(H)$. We note that line graph of any bipartite multigraph is a claw-free perfect graph. Moreover, the family consisting of the line graph of bipartite multigraphs forms a main role in characterization of claw-free perfect graphs \cite{MB}. In the light of Theorems \ref{first} and \ref{path cov} we should study Hamiltonian paths and path covering number of $\overline{L(H)}$. For convenience we call any graph of the form $\overline{L(H)}$ a skew graph, where $H$ is any bipartite multigraph. As explained earlier in the previous paragraph, the skew graphs can be constructed and represented in the following manner. Let $T$ be any $m\times n$ table of boxes containing $m$ rows and $n$ columns. Each row (resp. column) contains $n$ (resp. $m$) boxes. Assume that each box contains a non-negative integer. Hence, assume that the box placed in $i$-th row and $j$-th column contains an integer $t_{ij}$. Now, remove the integer $t_{ij}$ from the box and replace $t_{ij}$ distinct (and independent) vertices in the box. Hence, each vertex in the graph is placed in a box of the table. Let $u, v$ be two arbitrary vertices in the graph. Put an edge between $u$ and $v$ if and only if their corresponding boxes are in different rows and also in different columns. In fact, there exists no edge whose end-vertices belong to a same row (or to a same column). In this sense, each edge is a skew edge. The graph displayed in Figure \ref{fig-skew} is an example of skew graph and corresponding table is depicted in the figure. It's easily seen that $\overline{L(H)}$ can be represented by the above-mentioned method. The converse is also valid, i.e. any graph constructed as above is isomorphic to $\overline{L(H)}$ for some bipartite multigraph $H$. In fact, if all entries in the table $T$ is $0$ or $1$, then the corresponding graph $H$ is a simple bipartite graph.

\noindent The skew graphs can be expressed in terms of the concept of product dimension which was defined by Lov\'asz, Ne\v{s}et\v{r}il and Pultr in \cite{LNP}. Let $G$ be any graph without loops or multiple edges. According to \cite{LNP} the dimension of $G$ is the least $t$ for which there is a one-to-one function $f$ assigning to each vertex $v$ in $G$ a sequence $f(v)(1),f(v)(2), \ldots, f(v)(t)$ of real numbers so that $v$ is adjacent to $u$ in $G$ if and only if $f(v)(i)\not= f(u)(i)$ for $i=1,2, \ldots,t$. Observe that if in the definition of product dimension we allow general functions (and not necessarily one-to-one functions) then the collection of graphs with product dimension at most $2$ is equal to the family of skew graphs. The skew graphs are also appeared in \cite{HMT} (under the name of tabular graphs), where an application of these graphs in constructing graphs with high chromatic sum is shown.

\begin{figure}[ht]
\centering \unitlength=0.55mm \special{em:linewidth 0.4pt}
\linethickness{0.7pt}
\begin{picture}(45.00,120.00)
\put(25.00,70.00){\circle*{3.00}}
\put(25.00,64.00){$w$}
\put(25.00,88.00){\circle*{3.00}}
\put(25.00,106.00){\circle*{3.00}}
\put(25.00,109.00){$t$}
\put(5.00,60.00){\circle*{3.00}}
\put(45.00,60.00){\circle*{3.00}}
\put(-13.00,50.00){\circle*{3.00}}
\put(-19.00,50.00){$u$}
\put(63.00,50.00){\circle*{3.00}}
\put(65.00,50.00){$v$}
\emline{25.00}{70.00}{1}{25.00}{88.00}{2}
\emline{25.00}{70.00}{1}{5.00}{60.00}{2}
\emline{25.00}{70.00}{1}{45.00}{60.00}{2}
\emline{25.00}{106.00}{1}{25.00}{88.00}{2}
\emline{45.00}{60.00}{1}{25.00}{88.00}{2}
\emline{5.00}{60.00}{1}{25.00}{88.00}{2}
\emline{5.00}{60.00}{1}{45.00}{60.00}{2}
\emline{5.00}{60.00}{1}{-13.00}{50.00}{2}
\emline{63.00}{50.00}{1}{45.00}{60.00}{2}
\end{picture}\vspace{-2.5cm}
\caption{$u, v, w, t$ are tough vertices but no Hamiltonian path starts from any of them}\label{counter}
\end{figure}
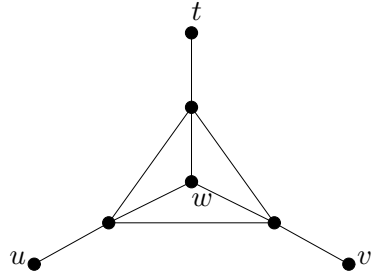

\noindent The relations between the concept of toughness and Hamiltonian cycles in graphs have been the research subject of many papers (e.g. \cite{DKS,LBZ,LWYZ,R,W}). In the sequel a graph $G$ is called {\it pathwise tough}
if $c(G\setminus S)\leq |S|+1$, for any $S\subseteq V(G)$,
where $c(G\setminus S)$ denotes the number of connected components of $G\setminus S$. Also, a vertex $v$ of $G$ is called {\it tough vertex} in $G$ if $(i)$ for any subset $S \subseteq (V(G)\setminus \{v\})$ we have $c(G\setminus \{v\} \setminus S) \leq |S|+1$ and $(ii)$ for any subset $S \subseteq (V(G)\setminus \{v\})$ such that $N(v)\subseteq S$ we have $c(G\setminus \{v\} \setminus S) \leq |S|$, where $N(v)$ is the neighborhood set of $v$ in $G$. It can be easily observed that if there exists a Hamiltonian path starting from a vertex $v$ of $G$ then $v$ is a tough vertex in $G$. But the converse is not true and it is easy to construct many counterexamples. In Figure \ref{counter} a graph is depicted in which the vertices $u,v,w,t$ are tough vertices but there exists no Hamiltonian path starting from any of them. Despite this situation, we prove in the rest of the paper that for any skew graph $G$, if $v$ is a tough vertex in $G$ then there exists a Hamiltonian path which starts from $v$.

\begin{prop}

\noindent $(i)$ If a graph $G$ contains a tough vertex then $G$ is pathwise tough.

\noindent $(ii)$ A vertex $v$ in a graph $G$ is tough vertex if and only if $G\setminus \{v\}$ is pathwise tough and for any $S\subseteq V(G)\setminus \{v\}$ with $N(v)\subseteq S$ we have $c(G\setminus \{v\} \setminus S)\leq |S|$.\label{tough-pathwise}
\end{prop}

\noindent \begin{proof}
 To prove part $(i)$, let $S$ be any subset of vertices in $G$. If $v\in S$ then $c(G\setminus S)=c(G\setminus \{v\} \setminus (S\setminus \{v\}))\leq |S|$. Now suppose that $v\not\in S$. If $N(v)\nsubseteq S$ then $\{v\}$ is not a connected component in $G\setminus S$ and hence $c(G\setminus S)\leq c(G\setminus (S\cup \{v\}))\leq |S|+1$. If $N(v)\subseteq S$ then $\{v\}$ is a connected component in $G\setminus S$ and $c(G\setminus S)\leq c(G\setminus\{v\} \setminus S)+1\leq |S|+1$. This completes the proof of part $(i)$. Proof of part $(ii)$ is clear.
\end{proof}

\noindent Let $G$ be any skew graph and $T$ its corresponding table. By a row $R$ (column $C$) in $T$ we mean the set of all boxes in $T$ which have a same row say $i$ (resp. a same column say $j$). Denote by $V(R)$ the set of all vertices placed in the boxes of $R$. Define $V(C)$ similarly. Note that each arbitrary row and column in the table intersects in exactly one box. We are going to obtain a necessary and sufficient condition in order that a given skew graph contains Hamiltonian path. Before we proceed we need to introduce four very spacial cases of tabular patterns. The patterns are depicted from left to right in Figure \ref{fig-pat} and are denoted symbolically by $\boxplus$, ${\mathcal{R}}$, ${\mathcal{C}}$ and ${\mathcal{R}}$$\cup$${\mathcal{C}}$, respectively. In all of the tables, each box containing $\ast$ means that there are some vertex or vertices in that box. Also any empty box means that the box contains no vertex. The first (from the left) is the $\boxplus$ pattern in which there are exactly four non-empty boxes in rectangular form as depicted in the figure. The pattern ${\mathcal{R}}$ (resp. ${\mathcal{C}}$) displays a typical table in which the vertices of the graph are placed only in one non-empty row (resp. column). In the unique non-empty row of ${\mathcal{R}}$, it is not necessary that all boxes are non-empty. The same is true for ${\mathcal{C}}$. The pattern ${\mathcal{R}}$$\cup$${\mathcal{C}}$ displays a typical skew graph in which all vertices are placed in the union of exactly one row say $R$ and one column say $C$. We exclude the case $V(R)\subseteq V(C)$ and also $V(C)\subseteq V(R)$ since otherwise the pattern reduces to either ${\mathcal{C}}$ or ${\mathcal{R}}$, respectively. The skew graph corresponding to any ${\mathcal{R}}$$\cup$${\mathcal{C}}$ pattern have $(V(R)\cap V(C))+1$ connected components. Here, we have used that $V(R)\setminus V(C)$ and $V(C)\setminus V(R)$ are both non-empty. The other patterns are also disconnected. We have the following theorem concerning disconnected skew graphs.

\begin{figure}
\begin{tabular}{|c|c|c|c|c|}
\hline&&&&\\[-1.30eM]
&&&&\\[.3eM]
\hline
~~&*&~~&*&~~\\[.3eM]
\hline
&&&&\\[.3eM]
\hline
&*&&*&\\[.3eM]
\hline
&&&&\\[.3eM]
\hline
\end{tabular}
\vspace{1cm}
\hspace{0.6cm}
\begin{tabular}{|c|c|c|c|c|}
\hline&&&&\\[-1eM]
&&&&\\[.3eM]
\hline
*&*&*&*&*\\[.3eM]
\hline
&&&&\\[.3eM]
\hline
&&&&\\[.3eM]
\hline
&&&&\\[.3eM]
\hline
\end{tabular}
\hspace{0.6cm}
\begin{tabular}{|c|c|c|c|c|}
\hline&&&&\\[-1.30eM]
&&*&~&\\[.3eM]
\hline
&&*&~&\\[.3eM]
\hline
&&*&~&\\[.3eM]
\hline
&&*&~&\\[.3eM]
\hline
~~&~~&*&~~&~~\\[.3eM]
\hline
\end{tabular}
\hspace{0.6cm}
\begin{tabular}{|c|c|c|c|c|}
\hline&&&&\\[-1eM]
&&*&&\\[.3eM]
\hline
&&*&&\\[.3eM]
\hline
*&*&*&*&*\\[.3eM]
\hline
&&*&&\\[.3eM]
\hline
&&*&&\\[.3eM]
\hline
\end{tabular}
\vspace{-0.5cm}
\caption{From left to right the patterns $\boxplus$, ${\mathcal{R}}$, ${\mathcal{C}}$, ${\mathcal{R}}$$\cup$~${\mathcal{C}}$}
\label{fig-pat}
\end{figure}

\begin{thm}
Let $G$ be a disconnected skew graph. Then tabular presentation of $G$ has one of the four patterns depicted in Figure \ref{fig-pat}.\label{patterns}
\end{thm}

\noindent \begin{proof}
Let the tabular presentation of $G$ be represented by the table $T$. An arbitrary box of $T$ is specified by a pair $(i,j)$, where $i$ (resp. $j$) is its row (resp. column). Also the entry of $T$ corresponding to $(i,j)$ denotes the number of vertices of $G$ whose location is that box. In case that this entry is non-zero we say that the box is non-empty. We claim that there do not exist three boxes in $T$ say $(i_1,j_1)$, $(i_2,j_2)$, $(i_3,j_3)$ such that $i_1, i_2, i_3$ are all distinct and $j_1, j_2, j_3$ are all distinct.
To prove the claim assume on the contrary that there exist boxes say $(i_1,j_1)$, $(i_2,j_2)$ and $(i_3,j_3)$ in $T$ such that $i_1, i_2, i_3$ are all distinct and $j_1, j_2, j_3$ are all distinct. We prove that $G$ is connected. First note that the subgraph of $G$ induced by the vertices of these three boxes is connected. Because no two boxes are in a same row or column. Now let $(i_0,j_0)$ be any arbitrary box in $T$ other than the former three boxes. It's clear that there exists $t\in \{1, 2, 3\}$ such that $i_0\not= i_t$ and $j_0\not= j_t$. It follows that all vertices corresponding to the box $(i_0,j_0)$ are adjacent to all vertices in the box $(i_t,j_t)$. It implies that the whole graph $G$ is connected. This contradiction proves the claim. Therefore there are at most two non-empty non-collinear boxes in $T$. Hence the non-empty boxes of $T$ can be covered by either two rows or two columns and or one row and one column. The third case is the same as the pattern ${\mathcal{R}}$$\cup$~${\mathcal{C}}$. If the boxes are covered by two rows then the pattern reduces either to $\boxplus$ or ${\mathcal{R}}$ or otherwise to a connected pattern. The argument for the case that non-empty boxes of $T$ are covered by two columns is exactly the same.
\end{proof}

\noindent Let $G$ be any skew graph and $T$ its corresponding table. Recall that by a row $R$ (column $C$) in $T$ we mean the set of all boxes in $T$ which have a same row say $i$ (resp. a same column say $j$). Also, $V(R)$ (resp. $V(C)$) is the set of all vertices placed in the boxes of $R$ (resp. $V(C)$). In the rest of paper the number of vertices in any graph $G$ is denoted by $|G|$.

\begin{thm}
Let $G$ be any connected skew graph and $T$ be its table. Then $G$ is pathwise tough if and only if the following conditions hold.

\noindent (i) For any row $R$ of $T$, $|G|\geq 2|V(R)|-1$.

\noindent (ii) For any column $C$ of $T$, $|G|\geq 2|V(C)|-1$.

\noindent (iii) For any row $R$ and column $C$ such that $V(R)\nsubseteq V(C)$ and $V(C)\nsubseteq V(R)$,\\
$|G| \geq |V(R)|+|V(C)|$.\label{skew-pathwise}
\end{thm}

\noindent \begin{proof} First assume that $G$ is not pathwise tough. Then there exists a subset $S$ of vertices such that $c(G\setminus S) \geq |S|+2$. Since $G$ is connected then $c(G)=1$. Hence $S$ is non-empty. If follows that $c(G\setminus S)\geq 3$. Note that $G\setminus S$ is a skew graph and since has at least $3$ connected components then by Theorem \ref{patterns} the table $T'$ corresponding to $G\setminus S$ obeys one of the patterns ${\mathcal{R}}$, ${\mathcal{C}}$ or ${\mathcal{R}}$$\cup$${\mathcal{C}}$ depicted in Figure \ref{fig-pat}. We check all cases.

\noindent {\bf Case 1.} Let the pattern of $T'$ be the pattern ${\mathcal{R}}$ and let its non-empty row be $R$.
Note that in this case $c(G\setminus S)=|V(G\setminus S)|=|V(R)|$. We have
$$|G|=|S|+|V(G\setminus S)|=|S|+c(G\setminus S)\geq 2|S|+2.$$
\noindent Now, $|G|=|S|+|V(R)|$ and $|G| \geq 2|S|+2$ imply $|G|\leq 2|V(R)|-2$ which contradicts the condition $(i)$.

\noindent {\bf Case 2.} Let the pattern of $T'$ be the pattern ${\mathcal{C}}$ and let its non-empty column be $C$. Using an argument similar to the case 1 we obtain $|G|\leq 2|V(C)|-2$, a contradiction with the condition $(ii)$.

\noindent {\bf Case 3.} In this case the pattern of $T'$ is pattern ${\mathcal{R}}$$\cup$${\mathcal{C}}$ of Figure \ref{fig-pat}. Let $R$ and $C$ be its non-empty row and column, respectively. Assume that there are exactly $p$ vertices in the intersection of $R$ and $C$. Note that in this case
$c(G\setminus S)=p+1$ and $|G|=|S|+|V(R)|+|V(C)|-p$. We have also $c(G\setminus S) \geq |S|+2$. Combining these relations we obtain $|G|\leq |V(R)|+|V(C)|-1$. But by the condition $(iii)$ we have $|G| \geq |V(R)|+|V(C)|$. This contradiction completes the first part of the proof.

\noindent Assume now that at least one of the conditions $(i)$, $(ii)$ or $(iii)$ does not hold. We prove that $G$ is not pathwise tough. We only consider the condition $(iii)$ because the proof for the other cases is similar. Let $R$ (resp. $C$) be an arbitrary row (resp. column) of $T$ with $V(R)\nsubseteq V(C)$ and $V(C)\nsubseteq V(R)$ such that $|G| \leq |V(R)|+|V(C)| -1$. Set $S= V(G)\setminus (V(R)\cup V(C))$. Assume that $|V(R)\cap |V(C)|=p$. We have
$$|G|=|S|+|V(R)\cup V(C)|=|S|+|V(R)|+|V(C)|-p \leq |V(R)|+|V(C)| -1.$$
\noindent It turns out that $|S|+1\leq p=c(G\setminus S)-1$. That is $c(G\setminus S) \geq |S|+2$, in other words $G$ is not pathwise tough.
\end{proof}

\noindent Let $v$ be any vertex in $G$ which is placed in row say $R_v$. Let $R$ be any row in the table of $G$ and $V(R)$ the set of vertices of $G$ which are in row $R$. If $R\not=R_v$ then $V(R)$ is the same as the set of vertices of $G\setminus \{v\}$ which are in row $R$. But if $R=R_v$ then the set of vertices in $G\setminus \{v\}$ which are in row $R$ is $V(R_v)\setminus \{v\}$. Also $N(v)\subseteq V(G)\setminus V(R)$ if and only if $R=R_v$. Now using Proposition \ref{tough-pathwise} and Theorem \ref{skew-pathwise} we obtain the following.

\begin{cor}
Let $G$ be any connected skew graph and $v$ any vertex of $G$ such that the tabular pattern of $G\setminus \{v\}$ is not the $\boxplus$ pattern. Let $T$ be the table of $G$. Assume that $v$ is located in row $R_v$ and column $C_v$ of $T$. Then $v$ is tough in $G$ if and only if the following conditions hold.
\begin{itemize}
\item{The $R_v$ condition: $|G|\geq 2|V(R_v)|-1$.}
\item{The $C_v$ condition: $|G|\geq 2|V(R_v)|-1$.}
\item{The $(R_v, C_v)$ condition: $|G|\geq |V(R_v)|+|V(C_v)|$.}
\item{The row condition: for each row $R\not=R_v$, $|G|\geq 2|V(R)|$.}
\item{The column condition: for each column $C\not=C_v$, $|G|\geq 2|V(C)|$.}
\item{The (row, column) condition: for each row $R$ and column $C$ with $(R,C)\not= (R_v,C_v)$,\\ $|G|\geq |V(R)|+|V(C)|+1$.}
\end{itemize}
\label{tough}
\end{cor}

\begin{prop}
Let $G$ be a skew graph and $v$ be a tough vertex in $G$. Let also $u$ be any neighbor of $v$ such that the table corresponding to the graph $G\setminus \{v,u\}$ has pattern $\boxplus$ in Figure \ref{fig-pat}. Then there exists a neighbor $w$ of $v$ in $G$ such that $w$ is tough in the graph $G\setminus \{v\}$. More strongly, $G$ contains a Hamiltonian path starting from $v$.
\label{square-case}
\end{prop}

\noindent \begin{proof}
It is enough to prove that $G$ contains a Hamiltonian path starting from $v$.
By the assumption, the tabular pattern of $G\setminus \{v,u\}$ is the $\boxplus$ pattern. Recall that in pattern $\boxplus$ there exist exactly four non-empty boxes. For the sake of simplicity, we denote these boxes by the set of vertices $A, B, C$ and $D$ in Figure \ref{fig-square-case}. Considering this $\boxplus$ pattern for $G\setminus \{v,u\}$, there are a few possibilities for the positions of $v$ and $u$ to be added in the $\boxplus$ pattern. We have searched for all possible cases and the result is the tables illustrated in Figure \ref{fig-square-case}. Note that $u$ and $v$ are not in a same row or column, because they are adjacent. Note also that it is possible that $v$ lies in one of the boxes corresponding to $A, B, C, D$. This is true for $u$ too, but this can not happen for $v$ and $u$ simultaneously. In other words, $\{v,u\}\nsubseteq A\cup B\cup C\cup D$. Because otherwise, the table of whole graph $G$ would be a $\boxplus$ pattern, i.e. a disconnected graph. Write $|A|=a, |B|=b, |C|=c, |D|=d$. We begin the argument by investigating the top-left pattern in Figure \ref{fig-square-case}.

\begin{figure}[t]
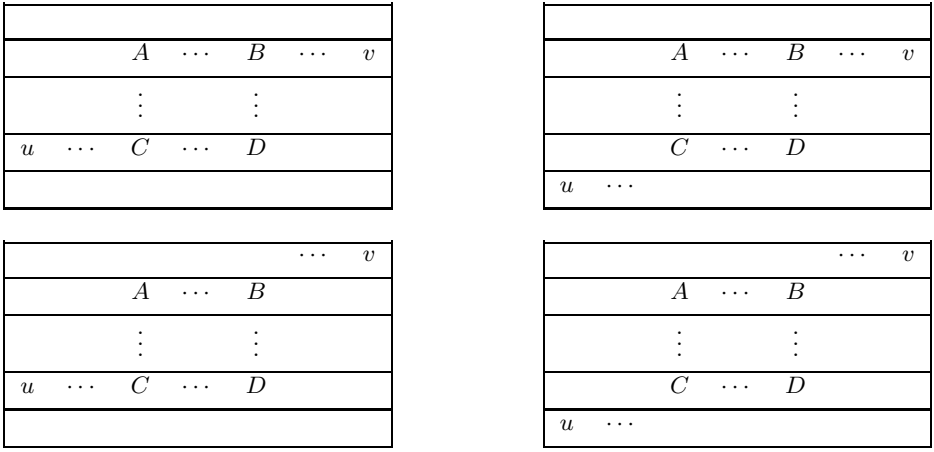

	\[
	{\footnotesize
		\hspace{0.6cm}
		\begin{tabular}{|ccccccc|}
		\hline&&&&&&\\[-1.30eM]
		&&&&&&\\[.3eM]
		\hline
		&~~&$A$&$\cdots$&$B$&$\cdots$&$v$\\[.3eM]
		\hline
		&&$\vdots$&&$\vdots$&&\\[.3eM]
		\hline
		$u$&$\cdots$&$C$&$\cdots$&$D$&&\\[.3eM]
		\hline
		&&&&&&\\[.3eM]
		\hline
		\end{tabular}
		\hspace*{2cm}
		\begin{tabular}{|ccccccc|}
		\hline&&&&&&\\[-1.30eM]
		&&&&&&\\[.3eM]
		\hline
		&~~&$A$&$\cdots$&$B$&$\cdots$&$v$\\[.3eM]
		\hline
		&&$\vdots$&&$\vdots$&&\\[.3eM]
		\hline
		&&$C$&$\cdots$&$D$&&\\[.3eM]
		\hline
		$u$&$\cdots$&&&&&\\[.3eM]
		\hline
		\end{tabular}}
	\]
	\vspace*{0.15cm}
	\[
	{\footnotesize
		\vspace{1cm}
		\hspace{0.6cm}
		\begin{tabular}{|ccccccc|}
		\hline&&&&&&\\[-1.30eM]
		&&&&&$\cdots$&$v$\\[.3eM]
		\hline
		&~~&$A$&$\cdots$&$B$&&\\[.3eM]
		\hline
		&&$\vdots$&&$\vdots$&&\\[.3eM]
		\hline
		$u$&$\cdots$&$C$&$\cdots$&$D$&&\\[.3eM]
		\hline
		&&&&&&\\[.3eM]
		\hline
		\end{tabular}
		\hspace{2cm}
		\begin{tabular}{|ccccccc|}
		\hline&&&&&&\\[-1.30eM]
		&&&&&$\cdots$&$v$\\[.3eM]
		\hline
		&~~&$A$&$\cdots$&$B$&&\\[.3eM]
		\hline
		&&$\vdots$&&$\vdots$&&\\[.3eM]
		\hline
		&&$C$&$\cdots$&$D$&&\\[.3eM]
		\hline
		$u$&$\cdots$&&&&&\\[.3eM]
		\hline
		\end{tabular}}
	\]
	\vspace*{0.3cm}
	\caption{The patterns involving in the proof of Proposition \ref{square-case}}
	\label{fig-square-case}
\end{figure}

\noindent {\bf The top-left pattern.}

\noindent First note that in this pattern we have also included the case in which the vertex $u$ is in the column of $B$ and $D$ (or the column of $A$ and $C$). The proof for this case is completely similar to the following proof for the top-left pattern. Also since this is the first pattern we investigate then we explain all details. Since $v$ is a tough vertex in the graph, we obtain the following inequalities.
$$c+d+1 = c(G\setminus (\{v\}\cup A \cup B)) \leq |A \cup B|+1= a+b+1.$$
$$a+b+1 = c(G\setminus (\{u\}\cup C \cup D)) \leq |C \cup D|+2= c+d+2.$$
$$b+d = c(G\setminus (\{v,u\}\cup A \cup C)) \leq 1+|A \cup C|+1= a+c+2.$$
$$a+c = c(G\setminus (\{v,u\}\cup B \cup D)) \leq 1+|B \cup D|+1= b+d+2.$$
$$c+1 = c(G\setminus (\{v\} \cup B)) \leq b+1.$$
$$d+1 = c(G\setminus (\{v\} \cup A)) \leq a+1.$$
\noindent Combining these inequalities concerning $a, b, c, d$, imply that either $b=c$ or $b=c+1$ and also either $a=d$ or $a=d+1$. In other words the distribution of vertices in the pair $(A, D)$ is balanced. The same is true for the pair $(B,C)$. We comment that this strong balance-ness result is valid for the rest of tables in Figure \ref{fig-square-case}. The main point is that a balanced distribution of vertices in $A$ and $D$ (also in $B$ and $C$) guarantees the existence of required Hamiltonian path. We come back to the top-left case again. In order to prove the existence of desired Hamiltonian path we check the following cases.

\noindent Case 1: $b=c$. We know that either $a=d$ or $a=d+1$. In this case we start from $v$ and go to $C$ and complete a zigzag path (round trip) between $C$ and $B$. Since $|B|=|C|$ this zigzag-like path ends at a vertex in $B$. Then we go from $B$ to $u$ and then go to a vertex in $A$. From $A$ we perform a zigzag path between $A$ and $D$. Since $|a-d|\leq 1$ this process outputs a path which visits all vertices, i.e. a Hamiltonian path. The constructed path can be denoted symbolically by $vCBuAD$. This notation for Hamiltonian paths is very helpful and will be used for the other patterns.

\noindent Case 2: $b=c+1$. In this case we have $a=d$. Construct the path as follows. Start from $v$ and go to $D$ then a zigzag path between $D$ and $A$ which ends at a vertex in $A$ (since $a=d$). We continue from $A$ to $u$ and then $B$ and complete the path by a zigzag path between $B$ and $C$. Similar to the previous case the path is denoted by $vDAuBC$.

\noindent {\bf The down-left pattern.}

\noindent In this case similar to the previous case we obtain the related inequalities between $a, b, c, d$. The results are $c\leq b \leq c+1$ and $d\leq a \leq d+1$. Consider first the case $b=c$. In this case the desired path is $vCBuAD$. Note that the path is valid for both cases $a=d$ and $a=d+1$. Consider now the case $b=c+1$. The Hamiltonian path is $vBCBuAD$ (for $a=d$) and $vBCBuADA$ (for $a=d+1$).
	
\noindent {\bf The top-right pattern.}
	
\noindent In this case we first obtain the relations $|b-c|\leq 1$ and $|a-d|\leq 1$. It follows then that the only possibilities are $(1)$ $b=c, a=d$, $(2)$ $b=c, a=d+1$, $(3)$ $b=c+1, a=d$, $(4)$ $b=c+1, d=a+1$. Note that in this case $u$ is adjacent to all vertices in $A\cup B \cup C \cup D$ and it is easy to obtain a desired Hamiltonian path. For example in the possibility $(1)$ the path is $vBCuAD$.
	
\noindent {\bf The down-right pattern.}

\noindent In this case too we obtain $|b-c|\leq 1$ and $|a-d|\leq 1$. Both $v$ and $u$ are adjacent to all vertices in $A\cup B \cup C \cup D$. Since the distribution of vertices in $B$ and $C$ (also in $A$ and $D$) are balanced. A desired Hamiltonian path is easily constructed.
\end{proof}

\section{Skew graphs and Hamiltonian paths: part II}

\noindent Now we state a necessary and sufficient condition in
order to have a Hamiltonian path in skew graphs. Also an
algorithm is presented which finds a Hamiltonian path (if any) in skew graphs.

\begin{thm}
Let $G$ be a connected skew graph and $v$ be any vertex of $G$. Then $G$ has a Hamiltonian path which starts from $v$ if and only if $v$ is a tough vertex in $G$.
\label{hamilton}\end{thm}

\noindent \begin{proof} \noindent First note that if $G$ has a Hamiltonian path starting from the vertex $v$ then $v$ is tough vertex. We prove by induction on the order of graph that if $G$ is any graph
and $v$ any tough vertex in $G$ then there exists a Hamiltonian path in $G$ which begins from $v$. We intend to use extensively and repeatedly Corollary \ref{tough} which provides a characterization of tough vertices. But one small and special pattern is excluded in Corollary \ref{tough}. This point is exactly concerned with the initial steps of our induction procedure. In order to overcome this issue we need to prove the following. If $G$ is any skew graph and $v$ is any tough vertex in $G$ such that for some neighbor $u$ of $v$ the tabular pattern of $G\setminus \{v,u\}$ is the $\boxplus$ pattern then there exists a Hamiltonian path in $G$ which starts from $v$. But this is exactly Proposition \ref{square-case}. Therefore based on this fact the initial steps of induction assertion hold. Assume now that the induction assertion holds for all graphs satisfying the specified conditions and with smaller than $n$ vertices. Let $G$ be a graph on $n$ vertices satisfying the conditions and let $v$ be any tough vertex in $G$. We first make the following claim.

\noindent {\bf Claim 1.} There exists a neighbor $w$ of $v$ such that $w$ is tough vertex in $G'=G\setminus \{v\}$.

\noindent {\bf Proof Claim 1.} Let $u$ be any neighbor of $v$. If $u$ is tough in $G'$ then $u$ is the desired vertex. Otherwise, $u$ is not tough in $G'$. Note that the tabular pattern of $G\setminus \{v,u\}$ is not the pattern $\boxplus$ (this is the initial step of induction), so we may apply Corollary \ref{tough}. We consider all conditions specified in Corollary \ref{tough} for $u$ and $G'$. Let $R_u$ and $C_u$ be the row and column corresponding to the location of $u$ in the table $G$, respectively.
\noindent \newline (1) The condition for $R_u$ is violated. In this case we have $|G'|\leq 2|V(R_u)|-2$. So $|G| \leq 2|V(R_u)|-1$. But $v$ is tough in $G$ and $v$ is not in $R_u$ (because $v$ and $u$ are adjacent). Hence the row $R_u$ has the same set of vertices in $G$ and $G'$. So we have $|G|\geq 2|V(R_u)|$. This is contradiction.
\noindent \newline (2) The condition for $C_u$ is violated. The argument in this case is exactly the same as previous one.
\noindent \newline (3) The condition for $(R_u,C_u)$ is violated. We have $|G'|\leq |V(R_u)|+|V(C_u)|-1$. But from other side we should have $|G|\geq |V(R_u)|+|V(C_u)|$, because $v$ is neither in row $R_u$ nor in column $C_u$. The inequalities contradicts each other.

\noindent Now one of the row, the column and the (row, column) conditions does not hold for the vertex $u$ in $G'$. We investigate each one separately.

\noindent {\bf (Case $1$) The row condition is violated for $(G',u)$:}

\noindent In this case there exists a row $R$ in the table of $G'$ such that $u\not\in V(R)$ and $|G'|\leq 2|V(R)|-1$. It follows that $|G|\leq 2|V(R)|$. We show that $v$ is not in the row $R$. Since otherwise the vertex $v$ is added to $V(R)$ in table of $G$ and $|G|\geq 2(|V(R)|+1)-1$. The inequalities are contradictory. So $v$ is not in row $R$ and hence $|G|\geq 2|V(R)|$. So $|G|=2|V(R)|$. Let $C_v$ be the column of $v$ in the table of $G$. We claim that $V(R)\nsubseteq V(R)\cap V(C)$. Otherwise the whole $V(R)$ is contained in $V(C_v)$ and hence $|V(C_v)|\geq |V(R)|+1$. Write the condition $C_v$ for $G$ and obtain $|G|\geq 2(|V(R)|+1)-1=2|V(R)|+1$. This contradiction shows that there exists a neighbor of $v$ say $w$ in the row $R$. The situation is illustrated in Figure \ref{fig-case1}. If $w$ is tough in $G\setminus \{v\}$ then the claim is proved. Otherwise, we check all conditions of Corollary \ref{tough} for the vertex $w$ and graph $G'$. Note that the row $R_w$ corresponding to $w$ is the very row $R$.

\begin{figure}
\setlength{\unitlength}{1.8mm}
\hspace*{-1.5cm}
\begin{picture}(80,40)
\put(40,8){\line(1,0){24}}\put(40,32){\line(1,0){24}}
\put(40,14){\line(1,0){24}}\put(40,18){\line(1,0){24}}
\put(57,25.5){$v$}
\put(35,15.5){$R$}\put(51,16){$w$}\put(45,22){$u$}
\large
\put(40,8){\line(0,1){24}}\put(64,8){\line(0,1){24}}
\small\put(35,2){} \large
\end{picture}
\vspace*{-1cm}
\caption{The situation of Case 1 in proof of Theorem \ref{hamilton}}
\label{fig-case1}
\end{figure}
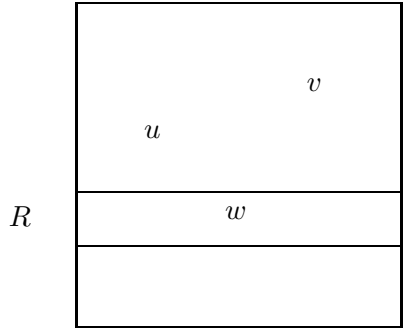

\noindent \newline (1) \underline{The condition for $R_w$ is violated.} We know $R_w=R$ and since $v\not\in R$ then the set of vertices in row $R$ is the same in both graphs $G$ and $G'$. From one side we have $|G|=2|V(R)|$ and from other side $|G\setminus \{v\}|\leq 2|V(R_w)|-2$. This is contradiction.

\noindent \newline (2) \underline{The condition for $C_w$ is violated.} We have $|G\setminus \{v\}|\leq 2|V(C_w)|-2$ and so $|G|\leq 2|V(C_w)|-1$. But sine $v$ is tough in $G$ then $|G|\geq 2|V(C_w)|$, a contradiction.

\noindent \newline (3) \underline{The condition for $(R_w,C_w)$ is violated.} In this case too we get a contradiction.

\noindent \newline (4) \underline{The row condition is violated for $w$ and $G'$.} In this case there exists a row $R'\not=R$ such that $|G\setminus \{v\}|\leq 2|V(R')|-1$. So $|G|\leq 2|V(R')|$. There are two possibilities. If $v$ is not in $R'$ then $|G|\geq 2|V(R')|$. It follows that $|G|=2|V(R)|=2|V(R')|$. Hence $V(G)=V(R)\cup V(R')$. This is a contradiction because $v$ is neither in $R$ nor $R'$. Assume that $v$ is in the row $R'$. In this case the set of vertices of $G$ in row $R$ has one vertex (i.e. the vertex $v$) more than $V(R')$. Since $v$ is tough it follows that $|G|\geq 2(|V(R')|+1)-1=2|V(R')|+1$, a contradiction.

\noindent \newline (5) \underline{The column condition is violated for $w$ and $G'$.}
In this subcase there exists a column $C$ consisting of $V(C)$ vertices in the table of $G\setminus \{v\}$ such that $|G\setminus \{v\}|\leq 2|V(C)|-1$. So $|G|\leq 2|V(C)|$. We know that $w$ is not in $C$. There are two possibilities. If $v$ is in the column $C$ then $|G|\geq 2(|V(C)|+1)-1=2|V(C)|+1$, a contradiction. If $v$ is not in column $C$ we have $|G|\geq 2|V(C)|$. Hence $|G|=2|V(C)|$. Now we apply the (row, column) condition for the vertex $v$ ad graph $G$ corresponding to $R$ and $C$. We obtain $|G|\geq |V(R)|+|V(C)|+1$. But $|G|=|V(R)|+|V(C)|$, a contradiction.

\noindent \newline (6) \underline{The (row, column) condition is violated for $w$ and $G'$.}
In this case there exist a row $R'$ and column $C$ in the table of $G\setminus \{v\}$ such that $w\not\in V(R')\cup V(C)$ and $|G\setminus \{v\}|\leq |V(R')|+|V(C)|$. Hence, $|G|\leq |V(R')|+|V(C)|+1$. We have $w\not\in V(R')\cup V(C)$ since similar to the previous arguments we obtain contradiction. Since $v$ is tough in $G$ we apply the (row, column) condition and obtain $|G|\geq |V(R')|+|V(C)|+1$, $|G|\geq 2|V(R')|$, $|G|\geq |V(R)|+|V(C)|+1$ and $|G|\geq 2|V(C)|$. These relations imply $|V(R)|=|V(R')|$ and so $V(G)=V(R)\cup V(R')$. There is now no row for the position of $v$ in $G$, a contradiction.

\noindent {\bf (Case $2$) The row condition is violated for $(G',u)$:}

\noindent In this case the argument is completely similar to that of the previous case.

\noindent {\bf (Case $3$) The (row, column) condition is violated for $(G',u)$:}

\noindent In this case we may assume that the row condition and column condition is not violated for $w$ and $G'$. Since the (row, column) condition is violated for $G'$ then there exists a row $R$ and column $C$ with $V(R)\nsubseteq V(C)$ and $V(C)\nsubseteq V(R)$ such that
$|G\setminus \{v\}|\leq |V(R)|+|V(C)|$. So $|G|\leq |V(R)|+|V(C)|+1$. Note that in this condition $u$ is placed neither in row $R$ nor in column $C$. A similar argument shows that position of $v$ can not be simultaneously in row $R$ and column $C$. The row condition holds for the vertex $u$ in $G'$, hence $|G\setminus \{v\}|\geq 2|V(R)|$. So $|G|\geq 2|V(R)|+1$ and similarly $|G|\geq 2|V(R)|+1$. Therefore $|G|=|V(R)|+|V(C)|+1$ and $|V(R)|=|V(C)|$. Assume that $v$ is not in the row $R$. Note that $V(R)\cap V(C)$ and $V(R)\setminus V(C)$ are both non-empty. Hence, there exists a neighbor of $v$ say $w$ in the row $R$. The situation is illustrated in Figure \ref{case3-pic}. If $w$ is tough in $G\setminus \{v\}$ then $w$ is the desired vertex and the claim is proved. Otherwise, one of the conditions six conditions is not satisfied for the graph $G\setminus \{v\} \setminus \{w\}$. We investigate each condition separately.

\begin{figure}
\hspace*{-2.5cm}
\setlength{\unitlength}{1.9mm}
\begin{picture}(90,40)
\put(40,8){\line(1,0){24}}
\put(40,12){\line(1,0){24}}
\put(40,16){\line(1,0){24}}
\put(40,24){\line(1,0){24}}
\put(40,28){\line(1,0){24}}
\put(40,32){\line(1,0){24}}
\put(55,25.5){$u$}
\put(55,33.5){$C'$}
\put(36,13.5){$R$}
\put(36,25.5){$R'$}
\put(45,13.5){$w$}
\put(45,33.5){$C$}
\put(61,29.5){$v$}
\large
\put(40,8){\line(0,1){24}}
\put(44,8){\line(0,1){24}}
\put(48,8){\line(0,1){24}}
\put(54,8){\line(0,1){24}}
\put(58,8){\line(0,1){24}}
\put(64,8){\line(0,1){24}}
\small\put(35,2){} \large
\end{picture}
\vspace*{-1cm}
\caption{The situation of Case 3 in proof of Theorem \ref{hamilton}}
\label{case3-pic}
\end{figure}

\noindent \newline (1) The condition for $R_w$ is violated. Note that $R_w=R$. We have $|G\setminus \{v\}|\leq 2|V(R_w)|-2$. So $|G|\leq 2|V(R_w)|-1$. But for the the vertex $v$ and graph $G$ and row $R$ ($v\not\in R$) we have $|G|\geq 2|V(R_w)|$, a contradiction.

\noindent The violation for $C_w$ and $(R_w,C_w)$ is carried out similarly. In the following we consider the violation of row condition, column condition and (row, column) condition as separate sub-cases.

\noindent {\bf Sub-case 1:} The row condition is violated for $w$ and $G'$. In this case there exists a row $R'$ (different from $R$) in the table of $G\setminus \{v\}$ such that $|G\setminus \{v\}|\leq 2|V(R')|-1$. So $|G|\leq 2|V(R')|$ and $w\not\in V(R')$. Concerning $R'$ and $v$ there are two possibilities. If $v\in R'$ then $|G|\geq 2(|V(R')+1)-1=2|V(R)|+1$, a contradiction. Therefore $v\not\in R'$ and so $|G|\geq 2|V(R')|$. It follows that $|G|=2|V(R')|$, i.e. an even integer. But we know that $|G|=2|V(R)|+1$. This contradiction rules out the possibility of this case.

\noindent {\bf Sub-case 2:} The prove for this case is similar to the previous one and we omit its details.

\noindent {\bf Sub-case 3:} In this case there exist a row $R'$ and column $C'$ in the table of $G\setminus \{v\}$ such that $|G\setminus \{v\}|\leq |V(R')|+|V(C')|$. The row $R'$ and $C'$ are indicated in Figure \ref{case3-pic}.
Since the conditions $(i)$ and $(ii)$ holds for $(G',w)$ (otherwise we go to either sub-case $1$ or sub-case $2$) then we use them for the row $R'$ and $C'$. We obtain $|G|\geq 2|V(R')|+1$ and $|G|\geq 2|V(C')|+1$ and hence $|G|=|V(R')|+|V(C')|+1$. Since $G\setminus \{v\}$ is pathwise tough then we have $|G|\geq |V(R)|+|V(C')|+1$ and $|G|\geq |V(R')|+|V(C)|+1$. It implies that $|V(R)|=|V(R')|$ and $|V(C)|=|V(C')|$. It follows that for some vertices $x$ and $y$, $V(G)=V(R)\cup V(R') \cup \{x\}$ and $V(G)=V(C)\cup V(C') \cup \{y\}$. Note that by their definition, $V(R)$ and $V(C)$ does not contain $v$ and $u$ because $V(R)$ is the set of vertices in $G\setminus \{v\} \setminus \{u\}$ whose position is in the row $R$. Similarly, $v\not\in V(R')\cup V(C')$. We conclude that $V(G)=V(R)\cup V(R') \cup \{v\}=V(C)\cup V(C') \cup \{v\}$. Hence, $u\in V(R')\cap V(C')$ and $w\in V(R)\cap V(C)$. More importantly, $V(G)\setminus \{v\} = (V(R)\cap V(C)) \cup (V(R)\cap V(C')) \cup (V(R')\cap V(C)) \cup (V(R')\cap V(C'))$. The situation is illustrated in Figure \ref{case3-pic}. It implies that the table of $G\setminus \{v\}$ obeys the pattern $\boxplus$. In other words, $G\setminus \{v\}$ is disconnected. This contradiction completes the proof of Claim $1$.

\noindent By Claim $1$, there exists a vertex $w$ such that $w$ is tough vertex in $G\setminus \{v\}$. The graph $G \setminus \{v\}$ satisfies the conditions of the theorem and has $n-1$ vertices. Hence by the induction it contains a Hamiltonian path which begins from $w$. This path together with the edge $vw$ forms a Hamiltonian path in $G$. This completes the proof of theorem.
\end{proof}

\noindent Let $G$ be a skew graph. Add a new vertex say $w$ to $G$ and connect $w$ to each vertex in $G$. The resulting graph $H$ is skew and $w$ is tough in $H$ if and only if $G$ is pathwise tough. Using this comment and applying Theorem \ref{hamilton} for $H$ and $w$ we obtain the following.

\begin{cor}
A skew graph contains a Hamiltonian path if and only if it is pathwise tough.
\label{hamilton-pathwise}
\end{cor}

\noindent The following results provide efficient algorithms to obtain and generate all Hamiltonian paths in skew graphs.

\begin{prop}
Assume that any skew graph is represented by its corresponding table of size say $p\times q$ consisting of non-negative integers. There exists an ${\mathcal{O}}(pqn)$ algorithm which decides whether a given skew graph $G$ and any tough vertex $v\in G$ contains a Hamiltonian path starting from $v$, where $n$ is the order of $G$. Moreover, the algorithm obtains one such Hamiltonian path in case that it exists.
\end{prop}

\noindent \begin{proof}
Let $G$ be any skew graph with a table of size $p\times q$ and $v$ be any vertex of $G$. By Theorem \ref{tough-pathwise}, $v$ is tough in $G$ if and only if $G\setminus \{v\}$ is pathwise tough. By Theorem \ref{skew-pathwise}, $G\setminus \{v\}$ is pathwise tough if and only if three conditions $(i)$ (for rows), $(ii)$ (for columns) and $(iii)$ (for pairs of rows and columns) hold. To check the validity of each single condition $(i)$, $(ii)$ and $(iii)$ needs at most $2pq$ time steps. Also there are $pq$ pairs of row and column. It follows that to check whether a given vertex is tough can be decided by ${\mathcal{O}}(pq)$ steps. At each step of constructing the desired Hamiltonian path in the proof of Theorem \ref{hamilton}, one vertex is removed from the underlying skew graph. We conclude that the Hamiltonian path is obtained with time complexity ${\mathcal{O}}(pqn)$.
\end{proof}

\noindent The proof of Theorem \ref{hamilton} provides more facilities about Hamiltonian paths in skew graphs. Let $G$ be a skew graph and $v$ a tough vertex in $G$. Existence of a Hamiltonian path from $v$ is guaranteed by Theorem \ref{hamilton}. By proof of Theorem \ref{hamilton} there exists a neighbor say $w$ of $v$ such that $w$ is tough vertex in $G\setminus \{v\}$. The vertex $w$ is obtained by an easy procedure in the proof. In fact all neighbors of $v$ which are tough in $G\setminus \{v\}$ are easily discovered. We pick an arbitrary neighbor of $v$ say $w$ which is tough. Then we obtain all neighbors of $w$ which are tough in $G\setminus \{v,w\}$. We continue this process and eventually reach at a Hamiltonian path in $G$. Note that by this method all Hamiltonian paths of $G$ beginning from $v$ are generated. The following result is yielded.

\begin{cor}
There exists an algorithm which generates all Hamiltonian paths which starts from $v$ for any given graph $G$ and tough vertex $v$ in $G$. Moreover, the generation of each individual path uses ${\mathcal{O}}(pqn)$ time steps, where $n$ is the order of $G$ and $p$ and $q$ are the sizes of the corresponding table of $G$.\label{generate}
\end{cor}

\noindent Recall that the path covering number $\pi(G)$ of $G$ is the smallest number of vertex disjoint paths needed to cover the vertices $G$. The following result of \cite{GMW} is reported in the introduction. If $\pi(\overline{G})\geq 2$ then
$\lambda_{2,1}(G)=|G|+\pi(\overline{G})-2$. The following result determines the path covering number of skew graphs. To present it we need a notation. Let $G$ be a skew graph and $T$ be its tabular presentation. For each row $R$ (resp. column $C$) of $T$, denote by $V(R)$ (resp. $V(C)$) the set of vertices whose location is in $R$ (resp. $C$). We associate a deficit parameter for $G$. For this reason we denote it tentatively by $DF(G)$ so that $G$ is pathwise tough if and only if $DF(G)=0$, i.e. in this case $G$ needs no more vertices for being pathwise tough. Let $R(T)$ and $C(T)$ be the sets of all rows and columns in the table $T$, respectively. In any term of the form $|V(R)|+|V(C)|$ in the following definition we assume that $V(R)\setminus V(C)$ and $V(R)\setminus V(C)$ are both non-empty.
$$DF(G)=\max \{2|V(R)|-1,~ 2|V(C)|-1,~ |V(R)|+|V(C)|:R\in R(T), C\in C(T)\}-|G|.$$

\begin{prop}
Let $G$ be a skew graph. Then $\pi(G)=DF(G)+1$.
\label{pathcov}
\end{prop}

\noindent \begin{proof}
Set $DF(G)=k$. Add a new row $R_{new}$ and column $C_{new}$ to $T$ and place $k$ independent vertices in the intersection of $R_{new}$ and $C_{new}$. Denote the resulting graph by $H$. Note that $H$ satisfies the conditions of Theorem \ref{skew-pathwise}, hence is pathwise tough. By Corollary \ref{hamilton-pathwise}, $H$ has a Hamiltonian path. This path introduces at most $k+1$ vertex disjoint paths in $G$ which cover all vertices of $G$. So $\pi(G)\leq DF(G)+1$.
An important point is that $k$ is the minimum number of vertices which should be added to $G$ in order to obtain a super-graph of $G$ satisfying the conditions of Theorem \ref{skew-pathwise}. Hence no Hamiltonian path in $H$ starts from any vertex in $V(R_{new})\cap V(C_{new})$. Based on this comment we prove
$\pi(G)\geq DF(G)+1$. Let $\{P_1, P_2, \ldots, P_t\}$ be the set of disjoint paths in a minimum path covering of $G$. Using these paths and the vertices of $S=V(R_{new})\cap V(C_{new})$ we obtain a Hamiltonian path in $H$. First traverse $P_1$ and then go to a vertex say $u_1$ in $S$. Then from $u_1$ to $P_2$ and traverse the whole $P_2$. Continue this procedure and obtain a Hamiltonian path $P$ in $H$. Note that based on the above comment, $P$ can not be finished at $S$. Therefore corresponding to the $k$ vertices of $S$ we require $k+1$ paths in $\{P_1, P_2, \ldots, P_t\}$, i.e. $t\geq k+1$.
\end{proof}

\noindent The following theorem reveals interesting facts concerning the path covering number of skew graphs. Denote by $\ell(G)$ the maximum length of any path in $G$. Let $G$ be a skew graph with $\pi(G)=t$ and $\{P_1, \ldots, P_t\}$ be any path covering for $G$. Then one of the paths is the path with maximum length in $G$ and the others are paths consisting of only one single vertex, i.e. of length zero. For example see the graph $G$ in Figure \ref{fig-skew}. We have $|G|=8$, $DF(G)=1$, $\pi(G)=2$ and $\ell(G)=6$. In fact $G$ contains a path on $7$ vertices.

\begin{thm}
In any skew graph $G$, $\ell(G)=|G|-DF(G)-1$.
\label{maxpath}
\end{thm}

\noindent \begin{proof}
Let $P$ be a path of length $\ell(G)$. Obviously $\pi(G)\leq 1+(|G|-\ell(G)-1)$. Hence $\ell(G)\leq |G|-\pi(G) =|G|-DF(G)-1$. The converse inequality is obtained by induction on $G$. For convenience set $M(G)=\max \{2|V(R)|-1, 2|V(C)|-1, |V(R)|+|V(C)|:R\in R(T), C\in C(T)\}$. It's easily observed that the maximum is taken for a unique row, or a unique column and or a unique pair of row and column. In each case there exists a vertex say $u$ such that $M(G')\leq M(G)-2$, where $G'=G\setminus \{u\}$. Note that $|G|-DF(G)=2|G|-M(G)\leq 2|G|-2-M(G')=2|G'|-M(G')=|G'|-DF(G')$. By applying the induction hypothesis for $G'$ we have $\ell(G')\geq |G'|-DF(G')-1\geq |G|-DF(G)-1$, as desired.
\end{proof}

\noindent Denote by $K_m\Box K_n$ the Cartesian product of $K_m$ and $K_n$. The $L(p,q)$-labeling including $\lambda$-labeling of $K_m\Box K_n$ have been the research subject of many papers (e.g. \cite{GMS,LLS}). The Cartesian product of other graphs has been studied in \cite{KY,SS1,SS}. In this regard we introduce the following mathematical item.

\begin{defin}
By an $m\times n$ $\lambda$-rectangle we mean any $m\times
n$ array with entries $1,2,\ldots mn$, such that in any row or
column any two entries differ at least two.
\end{defin}

\noindent An $m\times n$ $\lambda$-rectangle is equivalent to an $L(2,1)$-coloring of $K_m\Box K_n$. Using Corollary \ref{generate} we obtain the following.

\begin{thm}
For any $m\neq 1$ and except the case $m=n=2$, there is an
$m\times n$ $\lambda$-rectangle. Moreover, the algorithm of Corollary \ref{generate} generates all $\lambda$-rectangles.
\end{thm}

\section{Acknowledgment}

\noindent The author thanks the anonymous reviewers for their useful comments.



\begin{thebibliography}{1}
		
\bibitem{BKTV}
H.L. Bodlaender, T. Kloks, R.B. Tan, and J. van Leeuwen,
$\lambda$-coloring of graphs, Lecture Notes in Comput. Sci., 1770, Springer, Berlin, 2000, pp. 395--406.

\bibitem{BBFPW}
H. L. Bodlaender, H. Broersma, F. V. Fomin, A. V. Pyatkin, and G.
J. Woeginger, Radio labeling with pre-assigned frequencies, SIAM J. Optim. 15 (2004) 1--16.

\bibitem{DKS}
J. S. Deogun, D. Kratsch and G. Steiner, $1$-tough cocomparability
graph are Hamiltonian, Discrete Math. 170 (1997) 99--106.

\bibitem{FKK}
J. Fiala, T. Kloks, and J. Kratochvil, Fixed-parameter complexity
of $\lambda$-labellings, Discrete Appl. Math. 129 (2001) 59--72.

\bibitem{FKP1}
J. Fiala, J. Kratochvil, and A. Proskurowski, Distance
constrainded labellings of pre-colored trees, Lecture Notes in Comput. Sci., 2202, Springer, Berlin, 2001, pp. 285--292.

\bibitem{FNPS}
D. A. Fotakis, S. E. Nikoletseas, V. G. Papadopoulou, P. G. Spirakis, Hardness results and efficient approximations for frequency assignment problems: radio labelling and radio coloring, Comput. Inform. 20 (2001) 121--180.

\bibitem{FS}
D. A. Fotakis and P. G. Spirakis, A Hamiltonian approach to the assignment of non-reusable frequencies, Lecture Notes in Comput. Sci., 1530, Springer, Berlin, 1998, pp. 18--29.

\bibitem{GMS}
J.P. Georges, D.W. Mauro and M.I. Stein, Labeling products of
complete graphs with a condition at distance two, SIAM J. Discrete
Math. 14 (2001) 28--35.

\bibitem{GMW}
J.P. Georges, D.W. Mauro and M.A. Whittlesey, Relating path
coverings to vertex labelings with a condition at least two,
Discrete Math. 135 (1994) 103--111.

\bibitem{GY}
J.R. Griggs and R.K. Yeh, Labelling of graphs with a condition at
distance $2$, SIAM J. Discrete Math. 5 (1992) 586--595.

\bibitem{HMT}
H. Hajiabolhassan, M. L. Mehrabadi, R. Tusserkani, Tabular graphs
and chromatic sum, Discrete Math. 304 (2005) 11--22.

\bibitem{HCHY}
Y-Z. Huang, C-Y. Chiang, L-H. Huang, H-G. Yeh, On $L(2,1)$-labeling of generalized Petersen graphs, J. Combin. Optim. 24 (2012) 266--279.

\bibitem{KY}
D. Kuo and J.H. Yan, On $L(2,1)$-labeling of Cartesian products of
paths and cycles, Discrete Math. 283 (2004) 137--144.

\bibitem{LBZ}
B. Li, H.J. Broersma, S. Zhang, Forbidden subgraphs for Hamiltonicity of 1-tough graphs,
Discuss. Math. Graph Theory 36 (2016) 915--929.

\bibitem{LLS}
D. L\"{u}, W. Lin, Z. Song, Distance two labelings of Cartesian products of complete graphs,
Ars Combin. 104 (2012) 33--40.

\bibitem{LWYZ}
X. Li, B. Wei, Z. Yu, Y. Zhu, Hamilton cycles in $1$-tough triangle-free graphs,
Discrete Math. 254 (2002) 275--287.

\bibitem{LNP}
L. Lov\'asz, J. Ne\v{s}et\v{r}il, A. Pultr, On a product
dimension of graphs, J. Combin. Theory Ser. B 29 (1980) 47--67.

\bibitem{LZ}
C. Lu, Q. Zhou, Path covering number and $L(2,1)$-labeling number of graphs, Discrete Appl. Math. 161 (2013) 2062--2074.

\bibitem{MB}
F. Maffray, B.A. Reed, A description of claw-free perfect graphs, Journal of Combinatorial Theory, Series B 75 (1999) 134--156.

\bibitem{R}
L. Rao, A sufficient condition for 1-tough graphs to be Hamiltonian, in: Graph theory, combinatorics, and algorithms, Vol. 1, 2 (Kalamazoo, MI, 1992), 977--980, Wiley-Intersci. Publ., Wiley, New York, 1995.

\bibitem{SA}
U. Sarkar, A. Adhikari, On characterizing radio $k$-coloring problem by path covering problem,
Discrete Math. 338 (2015) 615--620.

\bibitem{SS1}
C. Schwarz, D. Sakai Troxell, $L(2, 1)$-labelings of products of two cycles, Disc. Appl. Math.154 (2006) 1522--1540.

\bibitem{SS}
Z. Shao, R. Solis-Oba, On some results for the $L(2,1)$-labeling on Cartesian sum graphs, Ars Combin. 124 (2016) 365--377.

\bibitem{SYZ}
Z. Shao, R.K. Yeh, D. Zhang, The $L(2,1)$-labeling on graphs and the frequency assignment problem, Applied Mathematics Letters 21 (2008) 37--41.

\bibitem{W}
B. Wei, Hamiltonian cycles in $1$-tough graphs, Graphs Combin. 12 (1996) 385--395.

\bibitem{Y}
R.K. Yeh, A survey on labeling graphs with a condition at distance two, Discrete Math. 306 (2006) 1217--1231.

\end{thebibliography}
\end{document}